\documentstyle[12pt]{article}
\long\def\notes#1{\ifinner
             {\tiny #1}
             \else
              \marginpar{\protect\tiny #1}%
              \fi}%
\begin{document}

\title{Use of an Hourglass Model in Neuronal Coding.}

\author{M. Cottrell, \thanks{ Postal address: SAMOS, Universit\'e Paris 1,
90, rue de 
Tolbiac, F-75634 Paris Cedex 13,  France}
\and T.S. Turova
\thanks{Postal address: University of Lund, Department of Mathematical
Statistics,
Box 118, S-221 00 Lund, Sweden}}
\date{}

\newcommand{\refs}[1]{{\rm (\ref{#1})}}
\renewcommand{\theequation}{\arabic{section}.\arabic{equation}}

\newtheorem{theo}{Theorem}
\newtheorem{defi}{Definition}
\newtheorem{pro}{Proposition}
\newtheorem{lem}{Lemma}
\newtheorem{rem}{Remark}
\newtheorem{cond}{Condition}
\newtheorem{cor}{Corollary}
\newtheorem{asu}{Assumption}

\maketitle

\bigskip
\begin{abstract}
We study a system of interacting
renewal processes which is a model for neuronal activity.
We show that the 
system possesses an exponentially large number (with respect to the 
number of neurons in the network) of  limiting configurations of the 
"firing neurons". These we call patterns.
Furthermore, under certain conditions of symmetry
we find an algorithm to control limiting patterns by means of the connection 
parameters.
\end{abstract}



\section{Introduction}
{\bf Definition of the model.}
We study here the {\it hourglass} model. This model 
is represented by a system of interacting
renewal processes,  numerated by the sites in
a finite  subset $\Lambda \subset {\bf
Z}^{\nu}$. Define
for any  $z \in \Lambda$ its neighbourhood:
\begin{equation}\label{dD}
 { D}(z) := \{ z \pm l_{k}, k=1, \ldots , K\} \cap \Lambda, 
\end{equation}
where $\{l_{k}, k=1, \ldots , K\} \subset {\bf
Z}^{\nu}$ is a fixed but arbitrary set 
of non-zero vectors.
Notice by our definition, that $z \not \in D(z)$, and for any $z,y \in  \Lambda$
\[z \in D(y) \mbox{ iff } y \in D(z).\]
Also let for any subset $J \subset \Lambda$
\[ { D}(J) := \{ j \pm l_{k}, k=1, \ldots , K, \ j \in J\}. \]

Define on some probability space $(\Omega, \Sigma, {\bf P})$ some
 random variables $X_{i}(0) $, $Y_i$, $\theta _{ij}$, $i \in {\Lambda}  $, $j \in D(i)$,
  to be
independent with densities having a certain number of
moments.
The variables $X_{i}(0)$
and $Y_i$, $i\in {\Lambda}$, are assumed to be positive and to represent the  
initial state and the self-characteristic of the $i$th neuron, 
respectively.
We shall call the
$\theta _{ij} $ the connection parameters. For any fixed  $i \in
{\Lambda}$, $j \in D(i)$, 
the distribution of $\theta _{ij} $ 
is assumed to be concentrated either on ${\bf R}_-$ or ${\bf R}_+$.
We will specify the
sign of  $\theta _{ij}$ in every example we treat below.
Assume also $\theta _{ij} \equiv 0$ for any $i \in {\Lambda}  $, $j \not \in D(i)$.
 Finally let $Y_i^{(n)}$, $Y_i^{(n,z)}$ and $\theta _{ij}^{(n)}$, $\theta
_{ij}^{(n,z)}$,
$n \geq 1$, $z \in \Lambda$, be independent
copies of the variables $Y_i$ and $\theta _{ij}$, respectively.

We shall define a Markov process $X(t)=(X_{i}(t), \, i \in
\Lambda )$, $t \geq 0$, with left-continuous
trajectories in ${\bf R}_{+}^{\Lambda}$ 
as follows. For all $i \in \Lambda$ and $t > 0$ define
\begin{equation}\label{MD}
X_i(t) = X_i(0)-t + \sum _{0 < t_n < t:\ X_i(t_n )=0 } \ Y_i^{(n)}
\end{equation}

\[+ \sum _{j \in D(i)} \ \ \sum _{0 < t_n < t:\ X_j(t_n )=0\ }
\left( (Y_i^{(n,j)}-X_i(t_n )) I \{ X_i(t_n ) \leq \theta _{ji}^{(n)} \} \right.\] 
\[\left. \hspace{3cm} - \theta _{ji}^{(n)} I \{ X_i(t_n ) > \theta _{ji}^{(n)} \} 
\right) \]

\[- \sum _{j \in D(D(i))\setminus \{i\} \cup  D(i) } \ 
\sum _{z \in D(i) \cap D(j):  \theta _{zi}<0} 
\sum _{
\begin{array}{rl}
0 < t_n < t: & X_j(t_n )=0, \\
 X_i(t_n )>\theta _{ji}^{(n)}, 
 \mbox{ and} & X_z(t_n ) \leq \theta _{jz}^{(n,j)} 
\end{array}
} 
\theta _{zi}^{(n,j)}.\]

\noindent
{\bf Biological interpretation.}
Any moment $t$ such that $X_z(t)=0$ for some $z \in \Lambda$ 
we call {\it the moment of firing} of the $z$th neuron.
Assume for a moment that $\theta _{ji}\equiv 0$ for any $i,j \in \Lambda
$. Then the last two summation terms
in the right-hand
side of  equation (\ref{MD})
 turn into zero,
meaning that $X_i(t)$ is merely a renewal process in this case.
Thus each
 component $X_{i}(t)$ represents {\it the duration of
time before the
next firing of the corresponding $i$th neuron assuming no
interaction takes place meanwhile}.

In the presence of non-zero interactions the dynamics 
 of  $X(t)=(X_{i}(t), $ $ i \in
\Lambda )$, $t \geq 0$, are
described
as follows.
As long as all the components of $X(t)$ are strictly
positive,
they decrease from the initial state $X(0)$
linearly in time with rate one until the first time
$t_{z_{1}}$
that one of the components reaches zero for some
$z_{1} \in {\Lambda }$: $
X_{z_{1}}(t_{z_1}) = 0$. We say that at this moment $t_{z_{1}}$,
the $z_{1}$th neuron fires and {\it sends impulse} $\theta _{z_1j} ^{(1)}$
to the $j$th neuron, if $j \in D(z_1)$. This means the
following. At the
 moment $t_{z_{1}}$, the trajectory $X_{z_{1}}(t_{z_{1}})$
jumps to a random value $Y_{z_1}^{(1)}$, i.e.
\begin{equation}\label{pt.3}
 X_{z_{1}}(t_{z_1}+) = Y_{z_1}^{(1)},
\end{equation}
which corresponds to the first summation
term in (\ref{MD}).
 At
the same moment every trajectory $X_{j}(t)$ with $j \in
D(z_{1})$ receives an increment $ \theta _{z_1j} ^{(1)} $
independent of the other processes, more precisely, according to the second
summation term in (\ref{MD})
\begin{equation}\label{0.1}
X_j(t_{z_1}+) =\left\{
\begin{array}{ll}
X_j(t_{z_1}) - \theta _{z_1j} ^{(1)}, & \mbox{ if }
X_j(t_{z_1}) - \theta _{z_1j} ^{(1)} >0, \\ \\

Y_j^{(1)}, & \mbox { otherwise. }
\end{array}
\right.
\end{equation}
Notice, that by (\ref{0.1}) the negative connections
$\theta _{z_1j}^{(1)}$  delay  the moment when $X_j(t)$ hits zero, 
i.e. the moment of firing of the $j$-th neuron.  Whereas positive connections
can shorten the time-interval until the next firing. That is why we shall call
negative  
connections {\it inhibitory} and the positive ones {\it excitatory}.
\medskip

{\bf Case 1.} Assume that all the 
interactions $\theta _{ij}$, $j \in {\it D}(i)$, are negative. 
Then definition (\ref{MD}) simply becomes 
\begin{equation}\label{ID}
X_i^{inh}(t) = X_i(0)-t + \sum _{0 < t_n < t:\ X_i(t_n )=0 } \ Y_i^{(n)}
+ \sum _{j \in D(i)} \ \ \sum _{0 < t_n < t:\ X_j(t_n )=0\ } \ |\theta _{ji}^{(n)}|,
\end{equation}
for $i \in \Lambda$ and $t \geq 0$.
We use a notation
$X^{inh}(t)$ for this particular case.
This model is due to Cottrell (1992) 
\cite{c} and it has been intensively studied
 (see an account of the previous results in Section 2 below).
\medskip

{\bf Case 2.} Assume, some of the interactions $\theta _{ij}$ are positive.
Notice that, only if $\theta _{z_1j_1}>0$,
it may happen that $X_{j_1}(t_{z_1}) - \theta _{z_1{j_1}} ^{(1)} \leq 0$. 
In this case  the trajectory $X_{j_1}(t)$ is reset instantaneously to an 
independent 
value $Y^{(1)}_j$ according to (\ref{0.1}). Furthermore, we say that
 the $j_1$th neuron also {\it fires } at the moment $t_{z_{1}}$ and 
changes instanteneously the trajectories of the neighbouring
 neurons according to the last term
in (\ref{MD}), which means the following.
 
Given $X_{z_1}(t_{z_1})=0$,
define the set of the firing at the same moment $t_{z_1}$ 
neurons:
\[
F_1 (z_1, t_{z_1}):=\{j_1 \in D(z_1): \ X_{j_1}(t_{z_1}) - \theta _{z_1j_1} ^{(1)} 
\leq 0 \}.
\]
Then at time ($t_{z_{1}}+$) the state of the system is defined as follows:
\begin{equation}\label{p.2}
X_{j}(t_{z_1}+)
\end{equation}
\[ =\left\{
\begin{array}{ll}
 Y_{j}^{(1)}, & \mbox{ if } \ j \in \{z_1\} \cup F_1, \\ \\
X_{j}(t_{z_1}) - \theta _{z_1j} ^{(1)}-
 \sum _i  \theta _{ij} ^{(1)},
& \mbox{ if } j \in  D(z_1) \setminus  F_1, \\ \\
X_{j}(t_{z_1}) - \sum _i \theta _{ij} ^{(1)},
 & \mbox{ if } j \in  D(F_1) \setminus \left( D(z_1) \cup  F_1 \cup \{z_1\} \right) ,
\end{array} \right.\]

\noindent
where $F_1=F_1 (z_1, t_{z_1})$,
the summation $\sum _i$ runs over the set 
$\{i \in F_1 \ : \  i \in D(j) \mbox{ and } \theta _{ij} <0 \}$,
and any summation over an empty set equals zero.
 The rest of the trajectories $X_{i}(t)$ with $i \not \in
D(F_1) \cup  F_1 \cup D(z_{1}) \cup \{z_{1}\} $  remain unchanged. After moment
$t_{z_{1}}$ the foregoing
dynamics are repeated.

To end this description recall the well-known fact from physiology that
a neuron does not react to the incoming impulses during
a certain 
period  right after 
it's own firing. This period  is called a refractory period.
 Observe that in our model a neuron receives an impulse only in the
moments when it does not fire itself. 
This reflects the property of the refractory
period.

This type of neural network has been proved  to be equivalent, in a
sense, to the "classical" neural model, which describes
the interacting membrane potentials (see \cite{t2}). 

\medskip
{\bf Plan of the paper.}
Our paper is organized as follows. In Section 2 we find
the critical values of the
parameters which separate the ergodic and transient cases.
 Also we provide some historical comments on our model in Section 2. 
In Section 3 we analyze a fully connected network with inhibitory
connections only, in which case we solve a problem related to the memory
capacity.

\section{Critical parameters.}
\setcounter{equation}{0}
\subsection{Results.}
Here we formulate our results on the critical values of the
parameters which separate the ergodic and transient cases
of the networks with excitatory and
inhibitory connections.

In order to eliminate boundary effects, we shall  assume here that 
$\Lambda = \{-N, \ldots,
N\}^{\nu}$ is a $\nu $-dimensional torus, i.e. we identify any
two points $(i_1, \ldots , i_{\nu})$ and $ (j_1, \ldots , j_{\nu})$ in 
${\bf Z}^{\nu} $
whenever
 $|i_m-j_m| \in \{0,2N\}$ for every $1 \leq m \leq \nu$. 
Hence, any point in $\Lambda$ has
the same
number of neighbours.
Further we 
assume that
$N>1$ and 
$\nu \in {\bf N}_+$ are fixed but arbitrary (although only the cases
$\nu =1,2,3$ are relevant for our context).
Define
for any $i, j \in \Lambda $
\[\|i-j\|_{\Lambda} := \sum_{k=1}^{\nu}|i_k-j_k|_{ 2N},\]
where 
$ |x-y|_{ 2N}:= \min \{|x-y|, \ |x-y-2N|, \ |x-y+2N|\} $ for any $x, y \in
{\bf Z}$.
Let $\Lambda _0$ be 
  the subset of the points in  $\Lambda$ such that:
\[
\begin{array}{rl}
i) & \mbox{the origin } (0, \ldots , 0) \in \Lambda _0, \\
ii) & \|x-y\|_{\Lambda} \neq 1 \mbox{ for any } x,y  \in \Lambda _0.\\
\end{array}
\]

Consider now the network $X(t)$ $=(X_i(t), i \in \Lambda )$ defined
in (\ref{MD}) with the following connection architecture.

\begin{asu}\label{a2} Let the neighbourhood $D(i)$ defined in (\ref{dD}) be
such that
\[D(i)=D_I(i) \cup D_E(i),\] 
where
\[ D_I(i)=\{j \in \Lambda: \ \|i-j\|_{\Lambda}=1\} \]
and
\begin{equation}\label{A2}
\begin{array}{ll}
 D_E(i) \subseteq \Lambda _0 & \mbox{ iff } i \in \Lambda _0, \\ 
 D_E(i) \subseteq \Lambda \setminus \Lambda _0 & \mbox{ iff } i
\in \Lambda \setminus \Lambda _0,
\end{array}
\end{equation}
with  $ | D_I(i) | =2\nu$ and $ |D_E(i)| = K_E$. The constants 
$0<K_E < N$ are fixed but arbitrary. 

\end{asu}
Then we set
\begin{equation}\label{1.1}
\theta _{ij}=\left\{
\begin{array}{ll}
-w_I  \, \eta _1^{ij}, & \mbox{ if } j \in D_I(i), \\ \\

w_E  \, \eta _2^{ij}, & \mbox{ if } j \in D_E(i),
\end{array}
\right.
\end{equation}
where $\eta _k^{ij}, i \neq j $ are independent copies of
 positive independent random variables $\eta _k $, respectively, $k=1,2,$
 with ${\bf E} \eta
_{k}=1$, and where  $w_I$ and
$w_E$ are  positive  parameters of the inhibitory and excitatory
connections, respectively. 
 
In words the condition (\ref{1.1}) means
that the nearest connections are inhibitory while the
more distant ones are excitatory.

Assume further that $X(0), \ X_i(0), i \in \Lambda $,
are {\it i.i.d.}, and  also 
$Y, \ Y_i, i \in \Lambda $,
are {\it i.i.d.} with ${\bf E}Y=1.$

\begin{asu}\label{a1} Assume, that
  the  distributions of $Y, \eta _1$, $
\eta _2$,  and  $X(0)$ have densities $g_0(u), g_1(u), g_2(u)$ and $ g_3(u)$,
respectively, which are positive differentiable functions, such that for 
some positive constants 
$a$ and $\alpha$ 
\begin{equation}\label{as1}
g_k(u) \leq a e^{-\alpha u} \ \ \mbox{ for all } u > 0 \mbox{ and } k=0,1,2,3.
\end{equation}
\end{asu}

We note that the
assumption of the exponential decay, though seemingly rather
restrictive,  arises naturally  from the physiology
 (see, for example \cite{t2} and the references
therein).

\begin{theo} \label{T1} 
Under Assumptions \ref{a2} and \ref{a1}
for any $w_E \geq 0$ there exists a positive $w_I^{cr}(w_E)$ such that the
system $X(t)$ with parameters (\ref{1.1}) is transient if
\begin{equation}\label{21.3}
w_I > w_I^{cr}(w_E),
\end{equation}
and ergodic if
\begin{equation}\label{21.33}
w_I < w_I^{cr}(w_E).
\end{equation}
\end{theo}

 We specify the critical function $w_I^{cr}(w_E)$ in (\ref{CR}) below.
Generally speaking, this function depends on $N$. However, its asymptotic
behaviour
when $w_E \rightarrow 0$ is uniform in $N$, as we establish
in the following theorem.

\begin{theo} \label{T2}
There exist positive constants $C$ and $C_0$ 
independent of $N$
such that
\begin{equation}\label{21.4}
\mid w_I^{cr}(w_E) \ - 
\ \left( \frac{1}{2 \nu }-\frac{K_E}{2 \nu } w_E \right) \mid  \leq Cw_E^2 
\end{equation}
for all $ 0 \leq w_E \leq C_0$.
\end{theo}

We postpone the proofs of these theorems to the next section.

So far only the case $w_E=0$ has been studied
analytically.  
Cottrell \cite{c} proved that when $w_E=0$, the network $X(t)=X^{inh}(t)$
 (see (\ref{ID}) ) is ergodic,  whenever
$w_I< 1/(2\nu )$,
which is, clearly  a particular case of (\ref{21.4}).
In \cite{c} 
sufficient conditions for convergence and transience  
were found for the finite network $X(t)$ with 
$\Lambda \subset {\bf Z}^2$ and 
 a specific structure of the  connections.
Further Piat \cite{p} (1994)
extended these results to a more general 
connection structure.
 Karpelevich {\it et al.} (1995) provided a
complete analysis of the evolution of inhibitory  networks when
the matrix of the
expected values of the connections defines 
a self-adjoint operator.

Paper \cite{cpr} presents simulated results for the model, in which  
both inhibitory and excitatory connections were incorporated. 
Here for the first time we 
 study analytically a rigorous mathematical model for such a network.
Our result explains the slope of the
diagram given on Fig. 4  \cite{cpr} in the neighbourhood of
the critical value $(0, 1/(2\nu) )$.
The most attractive feature of the hourglass network for neuromodelling
is that in the transient case 
the system splits as $t \uparrow \infty$
into two subsets. These are a subset of active (i.e. infinitely often firing) neurons 
and one of inactive neurons, which can be recognized as dark and white
areas, respectively, 
in the simulations, see e.g.
figures 3 and 6 of \cite{cpr}.  
A rigorous definition of the possible limiting
patterns of active and inactive neurons, called {\it traps} (see Section \ref{RS}
below) was 
given by Karpelevich, Malyshev and Rybko in \cite{kmr}. Thus
the set of all possible traps for a network should be
naturally thought of as a system of patterns,
which the network can hold, that is, memorize and recognize,
hopefully.
Recently \cite{mt} obtained a  result on phase transitions 
in the thermodynamical limit, which
establishes the large memory capacity of these networks. 

Further (in Section 3) besides a description of
the possible patterns, we solve an inverse problem. This is the problem of
how to determine the connections 
in order to get 
a system which stores  a given 
system of patterns.

\bigskip
\subsection{Proofs of Theorem  \protect \ref{T1} and Theorem  \protect \ref{T2}.}

\subsubsection{Preliminary definitions and results.}\label{RS}

We shall use in our proofs the results of \cite{kmr}.
Therefore firstly we introduce the terminology of \cite{kmr} for our model.

For any non-empty $W \subset \Lambda$ we call $ X^W(t) $ the
{\it restriction} of process  $X(t)$ on the set $W$. This is defined
as $X(t)$, but with initial conditions
\[ X_i^W(0) :=  \left\{
\begin{array}{ll}
X_i(0), & \mbox{ if } i \in W,\\
\infty,  & \mbox{ if } i \not \in W,
\end{array}
\right.
\]
in which case we assume that $X_i^W(t)=\infty$  if $i \not \in W$, for all $t>0$,
i.e. the nodes $\Lambda \setminus W$ are "deleted". 

Assume,
the process $ (X_i^W(t), i \in W) $ is ergodic for some non-empty 
$W \subset \Lambda$.
In this case we define for any $i \in W$ and $j \in D_E( i) \cap W$
 the following limiting frequencies:
\begin{equation}\label{f0}
\pi _{i}^{W, 0}:= \lim _{T \rightarrow \infty}
\frac{1}{T} \# \{0<t_n<T: \  X_i^{ W}(t_n)=0\},
\end{equation}
and 
\begin{equation}\label{fe}
\pi _{ij}^{W, e}:= \lim _{T \rightarrow \infty}
\frac{1}{T} \# \{0<t_n<T: \  X_j^{ W}(t_n)=0 \mbox{ and } X_i^{ W}(t_n)
- \theta^{(n)}_{ji} \leq 0\}.
\end{equation}
Thus $\pi _{i}^{W, 0}$ is the limiting frequency of 
firing of  the $i$th neuron
due to the hitting $0$ by the trajectory $X_i^{ W}(t)$, while 
$\pi _{ij}^{W, e}$ is the  limiting 
frequency 
of firing of the $i$th neuron due to an immediate excitatory impulse from the
$j$th neuron.  Then it is natural to call the {\it (total)
limiting frequency of firing} of
 the $i$th neuron the following sum of the defined above limits:
\begin{equation}\label{2.2}
\pi _{i}^{W}:= \pi _{i}^{W, 0} \ + \ \sum _{j \in D_E( i) 
\cap W} \ \pi _{ij}^{W, e}, \ \
i \in W.
\end{equation}
Next
we define the {\it second vector field} $v^{ W}=(v_j^{ W}, j \in 
\Lambda \setminus W)$  by
the following formula for its components:
\begin{equation}\label{2.1}
v_j^{W}=-1  - \sum_{i\in
D_I(j) \cap W} {\bf E} \theta_{ij} \ \pi _{i}^W \ - \ \sum_{i\in
D_E(j) \cap W} {\bf E} \theta_{ij} \ \pi _{i}^{W,0} , \ \ 
j \in 
\Lambda \setminus W . 
\end{equation}
Formally $v_j^{W}$ is the limiting mean drift of the 
$j$th component $X_j^W(t)$. Indeed, according to our definition (\ref{MD})
the $i$th neuron sends inhibitory impulses
to the corresponding neighbours with the limiting 
 frequency $\pi _{i}^{W}$, but it sends 
excitatory impulses with the limiting
frequency $\pi _{i}^{W, 0}$. 
For the details on the application of the
theory of the second vector field we refer to \cite{kmr} and the references
therein.

We shall call {\it a trap} for the process $X(t)$ 
any non-empty set $ M \subset \Lambda $ such that
the process
$(X_i^{ \Lambda  \setminus M}(t), i \in \Lambda  \setminus M)$ is ergodic,
while any coordinate of the second vector field $v^{\Lambda \setminus M}$
is positive, i.e.
\begin{equation}\label{2.21}
v_j^{\Lambda \setminus M} >0 \ \mbox{ for all } j \in M.
\end{equation}

For future reference let us
 rewrite now in our notations the criteria from \cite{kmr} on inductive ergodicity
and transience conditions.

\bigskip

\noindent
{\bf Theorem A} (See Theorem 2.1, \cite{kmr}.) {\bf I.} {\it Inductive ergodicity.}
If for any $W \subset \Lambda$ the process $(X_i^W(t), i \in W) $ is ergodic and
\[v_j^{ W} <0 \ \mbox{ for all } j \in \Lambda \setminus W\]
then the process $X(t)$ is ergodic.

\medskip

\noindent{\bf II.} {\it Sufficient transient conditions.}
If there exists a 
set $ M$ which is a trap for the process $X(t)$,
then the process $X(t)$ is transient.

\subsubsection{Subsystems with excitatory connections only.}

Consider now the restriction $X^{
\Lambda _0 }(t)$. Notice that
due to the definition of $\Lambda _0$ and Assumption \ref{a2}
the only connections
between the components of $(X_i ^{ \Lambda_0 }(t), i \in \Lambda_0)$ are
excitatory. 
Then by the definition of the process $X^{ \Lambda_0}(t)$ we have
for all $i \in \Lambda_0$:
\begin{equation}\label{1.101}
X_i^{ \Lambda_0}(t) = X_i(0)-t + \sum
_{0<s_n<t: \ X_i^{ \Lambda_0 }(s_n) = 0}Y^{(n)}
\end{equation}
\[ - \sum _{j \in D_E(i) } \ 
 \sum _{0<t_n<t: \ X_j^{ \Lambda_0 }(t_n) = 0 } \left( w_E\eta _{2}^{(n,j)} I\{ X_i^{
\Lambda_0 }(t_n) > w_E\eta _{2}^{(n,j)} \}
\right.\]
\[ \left. +( X_i^{ \Lambda_0 }(t_n)-Y ^{(n,j)}) I\{ X_i^{
\Lambda_0 }(t_n) \leq  w_E\eta _{2}^{(n,j)} \} \right), \]
where $\eta _{2}^{(n,j)}$, $n \geq 1,$ $j \in \Lambda _0$, 
are independent copies of the variable $\eta _{2}$.
Clearly, for any $i \in \Lambda_0$ the one-dimensional $X_i^{\{i\}}(t)$  is ergodic. 
Further for any $w_E>0$ and for any $W \subset \Lambda_0$ such that
 $(X_i^W(t), i \in W)$ is ergodic, 
the second vector field
defined in (\ref{2.1}) has negative coordinates only:
\[
v_j^{W}=-1 - \ \sum_{i\in
D_E(j) \cap W} w_E \ \pi _{i}^{W,0} , \ \ 
j \in 
\Lambda _0 \setminus W . 
\]
Hence,
by Theorem A  on inductive ergodicity 
$(X_i ^{ \Lambda_0 }(t), i \in \Lambda_0)$ is ergodic for any fixed $w_E>0$.
Taking into account translation invariance,  
we derive from (\ref{1.101}) and (\ref{2.2})
\begin{equation}\label{1.5}
\pi _i^{  \Lambda_0 } =:\pi ^+(w_E)=
\pi ^{+, 0}(w_E) \ + \ K_E \ \pi ^{+, e} (w_E)>0
\end{equation}
for any $i \in  \Lambda_0 $.
In particular,
\begin{equation}\label{1.9}
\pi ^+(0)=\pi ^{+, 0}(0) = 1/{\bf E}Y=1,
\end{equation}
 since in the case
$w_E=0$ we have a system of independent renewal processes
numerated by the sites of $ \Lambda_0$. 

Thus one can view the system $X(t)$ as a result of inhibitory 
interactions
between two (independent at initial moment) excitatory networks: $(X _i(t),$ $ \ i \in
\Lambda_0 )$ and
$(X_i (t), \ i \in \Lambda \setminus \Lambda_0)$.

\subsubsection{ Proof of Theorem  \protect \ref{T1}.}

We shall show that under condition 
(\ref{21.3}), with
\begin{equation}\label{CR}
w_I^{cr}(w_E)= \frac{1}{2 \nu \pi ^+(w_E)}
\end{equation}
the set $\Lambda _{0}$ is a trap for the system $X(t)$. Then the 
first statement on transience in our theorem follows from Theorem A.

As we have shown above, $(X^{\Lambda _{0}}(t), \ i \in \Lambda _{0})$ is ergodic.
Computing the second vector field  
with respect to formula (\ref{2.1}), and taking into account (\ref{1.1}) and
(\ref{1.5}) we get
\begin{equation}\label{1.15}
v_j^{ \Lambda_0  } = -1  + \sum_{i\in
D_I(j) }w_I \pi _{i}^{\Lambda_0 }=-1  + 2 \nu w_I \pi^+(w_E), \ \ 
j \in \Lambda \setminus \Lambda_0.
\end{equation}
Hence condition (\ref{21.3}) together with
(\ref{CR})
is equivalent to
\begin{equation}\label{1.16}
v_j^{ \Lambda_0  }>0, \ \ 
j \in \Lambda \setminus \Lambda_0,
\end{equation}
i.e. inequality in  (\ref{2.21}) is fulfilled for all $j \in \Lambda
\setminus \Lambda _{0}$, which implies that  $\Lambda _{0}$ is a trap. Hence, 
 our statement on transience follows by  Theorem A.

Assume now that condition (\ref{21.33}) with (\ref{CR}) holds.
Clearly, for any $i \in \Lambda$ one-dimensional $X_i^{\{i\}}(t)$  is ergodic. 
Next we will show that for any ergodic face $W$ the components of the second
vector field $v^{W}$ defined by (\ref{2.1})
are negative. It obviously follows from our definition of the model,
 that the limiting firing frequency is the highest
when the neuron has the excitatory connections only 
and the maximal number of them, i.e.
\[\pi^W_i \leq \pi _{j}^{\Lambda_0 }=\pi^+(w_E)\]
for any $W \in \Lambda$, any $i \in \Lambda \setminus W $ 
and $ j \in \Lambda \setminus \Lambda _{0}$. From here and (\ref{2.1})
we derive for any 
$i \in \Lambda \setminus W $ 
\[v_i^{W} \leq   -1  + 2 \nu w_I \pi^+(w_E)<0,
\]
where the last inequality is due to  (\ref{21.33}).  Hence,
using Theorem A (Theorem 2.1 \cite{kmr}) on inductive ergodicity, we readily derive our
statement on ergodicity. Theorem \ref{T1} is proved.

\subsubsection{Proof of Theorem  \protect \ref{T2}.}
Let $p_A(t,u_A)$, $u_A \in {\bf R_+}^{|A|}$, denote for any
finite $A \subseteq   \Lambda_0 $  the density of the
process  $X_A^{\Lambda_0}(t)=(X_i^{\Lambda_0}(t), i \in A)$. In the particular 
case $w_E=0$, we denote the corresponding density by $p_A^{0}(t,u_A)$.
Due to ergodicity 
the following limit
\begin{equation}\label{1.7}
\lim_{t \rightarrow \infty} p_A(t,u_A)=:p_A(u_A)
\end{equation}
exists for any $A \subseteq \Lambda_0 $. In particular, we have 
 for the limiting density
of an independent  renewal process:
\begin{equation}\label{1.700}
\lim_{t \rightarrow \infty} p^0_A(t,u_A)=p^0_A(u_A)=
\prod_{j \in A}p_j^0(u_j)
\end{equation}
where
\begin{equation}\label{1.70}
p_j^0(u) = \frac{\int_{u}^{\infty}g_0(v)dv}{{\bf E}Y}=:p^0(u).
\end{equation}

Further we will use the following lemma stating that
the limiting density in (\ref{1.7}) is a small perturbation of the function in 
(\ref{1.700})
when the parameter of interaction $w_E$ is sufficiently small.

\begin{lem} \label{L2} There exist positive and independent of $N$ constants 
$C, \beta $ and $c$ such that for any $0< w_E \leq c$ and for any 
subset      $A \subseteq \Lambda_0 $
\begin{equation}\label{2.L2}
| p_{A}(t, u_{A})- p_{A}^0(t, u_{A}) | \leq w_E C^{|A|}
e^{ - \beta \sum_{i \in A} u_i}, 
\end{equation}
for any $t >0 $ and $u_{A} \in {\bf R}_+^{A}$.
\end{lem}

{\bf Proof.} Recall that $g_0$ and $g_3$ are  the density functions
 of the variables $Y$ and $X(0)$, respectively.
Let $p(t, u_{\Lambda _0}, v_{\Lambda _0})$, $u_{\Lambda _0}, v_{\Lambda _0} 
\in {\bf R_+}^{\Lambda_0}$, denote the transition density of the 
process $(X_i^{\Lambda _0}(t), \ i \in \Lambda _0 )$. Then for any
 $A \subseteq   \Lambda_0 $ 
the density  $p_A(t,u_A)$, $u_A \in {\bf R_+}^{|A|}$, is given by the formula
\begin{equation}\label{red}
p_A(t,u_A) = \int _{{\bf R_+}^{|\Lambda_0 \setminus A|}} \
\int _{{\bf R_+}^{\Lambda_0}} \left( \prod_{z \in \Lambda _0} g_3(u_z)\right)
\ p(t, u_{\Lambda _0}, v_{\Lambda _0}) \ d\,u_{\Lambda _0} \ 
d\,v_{\Lambda_0 \setminus A}.  
\end{equation}

In the particular case $w_E=0$, we denote the transition density
by $p^{0}(t, u_{\Lambda _0}, v_{\Lambda _0})$. Clearly, in this case
\[p^{0}(t, u_{\Lambda _0}, v_{\Lambda _0})=\prod_{j \in \Lambda _0}
\, p^{0}(t, u_j, v_j)\]
due to the independence of the components. Notice, that the transition
 density of each of the
independent renewal processes is
\begin{equation}\label{ren}
p^{0}(t, u, v) \, = \, 
\end{equation}

\[ = \left\{ \begin{array}{ll}
\delta (u-t-v), & \mbox{if $0 \leq t < u$,} \\
 g_0(t-u+v)\, + \, \sum _{k=1} ^{\infty} \, \int _{0} ^{t-u}p_{S_{k}}(x)g_0(t-u+v-x) \, 
dx, & \mbox{if $0 \leq u \leq t$,}
\end{array}
                \right.\] 
for every $u,v \in {\bf R}_{+}$ and $t>0$, where $\delta (u \, - \, \cdot) $ 
is the Dirac distribution concentrated at $u$; $Y, Y^{k}, \, k \geq 1$,
are {\it i.i.d.};
$S_{k} := \sum _{l=1} ^{k} \, Y^{l}, \, k \geq 1$,
 and $p_{S_{k}}$ is the density of 
the distribution $S_{k}$. 

Using the results of Stone \cite{C.S.} one can show that under Assumption \ref{a1}
there exist positive constants $C'$ and $\alpha '$ such that
\begin{equation}\label{p1}
|p^{0}(t+t', 0, v) -p^{0}(t, 0, v)|\, \leq \, C'e^{- \alpha '(t+v)}
\end{equation}
for all $v \in {\bf R}_{+}$ and $t, t'>0$ (for the details see for example, \cite{t1}
Section 4). Combining (\ref{p1}) and (\ref{1.700}), we get:
\begin{equation}\label{ph2}
|p^{0}(v)- p^{0}(t, 0, v) |\, \leq \,  C'e^{- \alpha '(t+v)}.
\end{equation}
(For the general theory of piece-wise linear stochastic processes we refer to
the book \cite{D} by Davis, 1993.)

After these preliminaries the proof of 
Lemma \ref{L2} follows by standard techniques using
so-called cluster expansions (see, e.g.  \cite{t1} for a similar model
and \cite{mit} for a class of stochastic processes with local weak interactions)
as soon as the Kolmogorov equation (see (\ref{re.1}) below)
is written down for 
the transition density $p(t, u_{\Lambda _0}, v_{\Lambda _0})$.

Let ${\it M}({\bf R}_{+}^{\Lambda _0})$ denote a class of densities 
$\mu $  on ${\bf R}_{+}^{\Lambda_0}$ 
which possess all 
continuous partial derivatives. Define also
\[p(t,u_{\Lambda _0}, \Gamma)= \int_{\Gamma}p(t,u_{\Lambda _0},v_{\Lambda _0}) 
dv_{\Lambda _0}\]
for any Borel set $\Gamma \in {\bf R}_+^{\Lambda _0}$.
Denote here ${\it O}(i)=D_E(i)$ and $w_E =\epsilon $.
Then for our model we derive for any 
density $\mu \in {\it M}({\bf R}_{+}^{\Lambda _0})$  and for any 
Borel set $\Gamma \in {\bf R}_+^{\Lambda _0}$ 
\begin{equation}\label{re.1}
\frac{\partial }{{\partial } t } \int _{{\bf R}_{+}^{\Lambda _0}} \,  
\mu (u_{\Lambda _0}) \,\, 
p(t, u_{\Lambda _0}, \Gamma) \, d  u_{\Lambda _0} \, 
\end{equation}
\[= \, 
\, \,  - \, \int _{{\bf R}_{+}^{\Lambda _0}}\, \, 
\sum _{i \in \Lambda _0}  \, \,  \,
 \mu (u_{\Lambda _0}) \, \frac{\partial }{{\partial } u_{i} }  \, 
p(t, u_{\Lambda _0}, \Gamma)\, d  u_{\Lambda _0}  \,   \]
\[ \begin{array}{ll}
+ &
\int _{{\bf R}_{+}^{\Lambda _0}} \, \, 
 \mu (u_{\Lambda _0}) \, \sum _{i \in \Lambda _0}  \, \, 
\delta ( u_{i})\, \,
  \int _{{\bf R}_{+}^{1+2|{\it O(i)}|}} \, 
p_Y( w_{i}) \, \left( \prod _{j\in {\it O}(i)}p_Y(y_j) \, p_{\epsilon  \eta_2}(w_{j}) 
\right) \\ \\

& \times \, [\, p (t, u_{\Lambda _0} + U_{\Lambda _0}(u_{\Lambda _0}, 
w_{i \bigcup {\it O(i)}},y_{{\it O(i)}}), \Gamma) \, - \, 
p (t, u_{\Lambda _0}, \Gamma) \, ] \\ \\

& \times \,  d\, y_{{\it O(i)}} \  d \, w_{i \bigcup {\it O(i)}}  \,  \, 
d  u_{\Lambda _0},  
\end{array}\]
where $\delta (\cdot) $ is the Dirac measure concentrated at $0$,
 \[p_Y(u)=g_0(u), \ \ p_{\epsilon  \eta_2}(u)=\left\{ 
\begin{array}{ll}
g_2(u/{\epsilon })/{\epsilon },  & \mbox{ if } \epsilon  >0,\\
\delta(u),  & \mbox{ if } \epsilon  =0,
\end{array}
\right.
\ \ \ u \geq 0;\]
and the vector
$U_{\Lambda _0}(u_{\Lambda _0}, w_{i \bigcup {\it O(i)}},y_{{\it O(i)}})$ 
has  components:
\[ U_j(u_{\Lambda _0}, w_{i \bigcup {\it O(i)}},y_{{\it O(i)}})\]
\[=\left\{
\begin{array}{ll}
w_i, & \mbox{ if } j=i,\\
-w_jI\{w_j<u_j\} +(y_j-u_j)I\{w_j \geq u_j\}, & \mbox{ if } j\in {\it O(i)},\\
0, & \mbox{ otherwise. }
\end{array}
\right. \]

Further for any $w_i \geq 0$ let the vector ${\bar w_i} =({\bar w_i} \, _j, j \in
\Lambda _0)$ have the following components:
\[{\bar w_i}\, _j=\left\{
\begin{array}{ll}
w_i, & \mbox{ if } j=i,\\
0, & \mbox{ otherwise. }
\end{array}
\right.\]
Let us rewrite now formula (\ref{re.1}) in the following operator form:
\begin{equation}\label{re.2}
\frac{\partial }{{\partial } t } \int _{{\bf R}_{+}^{\Lambda _0}} \,  
\mu (u_{\Lambda _0}) \,\, 
p(t, u_{\Lambda _0}, \Gamma) \, d  u_{\Lambda _0} \, 
\end{equation}
\[=  \, \int _{{\bf R}_{+}^{\Lambda _0}}\, \, 
\mu (u_{\Lambda _0}) \, \sum _{i \in \Lambda _0} \,
 \left( 
  - \, \frac{\partial }{{\partial } u_{i} }  \, 
p(t, u_{\Lambda _0}, \Gamma)  \right. \]
\[ \left. + \, \delta ( u_{i})\, 
\int _{{\bf R}_{+}} \,  \, 
p_Y( w_{i})[ p (t, u_{\Lambda _0} + \bar{w_i}, \Gamma) -
p (t, u_{\Lambda _0}, \Gamma)]\, d  w_{i} \, \right)\,
 d u_{\Lambda _0}  \]

\[ \begin{array}{ll}
+ &
\int _{{\bf R}_{+}^{\Lambda _0}} \, \, 
 \mu (u_{\Lambda _0}) \, \sum _{i \in \Lambda _0}  \, \, 
\delta ( u_{i})\, \,
  \int _{{\bf R}_{+}^{1+2|{\it O(i)}|}} \, 
p_Y( w_{i}) \, \left( \prod _{j\in {\it O}(i)}p_Y(y_j) \, p_{\epsilon  \eta_2}(w_{j}) 
\right) \\ \\

& \times [\, p (t, u_{\Lambda _0} + U_{\Lambda _0}(u_{\Lambda _0}, 
w_{i \bigcup {\it O(i)}},y_{{\it O(i)}}), \Gamma) \, - \, 
p (t, u_{\Lambda _0} + \bar{w_i}, \Gamma) \, ] \\ \\

& \times \,  d\, y_{{\it O(i)}} \  d \, w_{i \bigcup {\it O(i)}}  \,  \, 
d \,  u_{\Lambda _0}  
\end{array}\]

\[=: \int _{{\bf R}_{+}^{\Lambda _0}} \, \, 
 \mu (u_{\Lambda _0}) \, \sum _{i \in \Lambda _0} (H^{i,\Lambda _0}_0
+H^{i,\Lambda _0}_1) \,
p (t, u_{\Lambda _0}, \Gamma) d \,  u_{\Lambda _0}. \]
Notice, that when $\epsilon =0$ the operator $H^{i,\Lambda _0}_1 \equiv 0$
for any $i$. In this case (\ref{re.2}) simply becomes 
\begin{equation}\label{re.3}
\frac{\partial }{{\partial } t } \int _{{\bf R}_{+}^{\Lambda _0}} \,  
\mu (u_{\Lambda _0}) \,\, 
p^0(t, u_{\Lambda _0}, \Gamma) \, d  u_{\Lambda _0} \, 
= \int _{{\bf R}_{+}^{\Lambda _0}} \, \, 
 \mu (u_{\Lambda _0}) \, \sum _{i \in \Lambda _0} \, H^{i,\Lambda _0}_1 \, 
p^0 (t, u_{\Lambda _0}, \Gamma)\, d \,  u_{\Lambda _0}
\end{equation}
i.e., the equation for the 
density of the process whose components are independed renewal processes.
Hence the formula (\ref{ren}) gives a solution to (\ref{re.3}):
\[p^0 (t, u_{\Lambda _0}, \Gamma)= \int _{ \Gamma} \, p^0 (t, u_{\Lambda _0}, 
v_{\Lambda _0}) \, d \,v_{\Lambda _0}. \]
Now we can find the solution to equation (\ref{re.2}) 
by solving the equivalent integral equation
(for the reference see also \cite{t1} and \cite{mit}):
\begin{equation}\label{re.4}
\int _{{\bf R}_{+}^{\Lambda_0}}\,
\, \mu (u_{\Lambda_0}) \, p(t, u_{\Lambda_0}, \Gamma) \, d u_{\Lambda_0}
\, = \, \int _{{\bf R}_{+}^{\Lambda_0}}\,
\, \mu (u_{\Lambda_0}) \, p^{0}(t, u_{\Lambda_0}, \Gamma) \, d u_{\Lambda_0} 
\end{equation}
\[\, + \, \int _{0}^{t} \int _{{\bf R}_{+}^{\Lambda_0}} 
\, \int _{{\bf R}_{+}^{\Lambda_0}}\,
\, \mu (u_{\Lambda_0}) \,
p^{0}(t-s, u_{\Lambda_0}, u_{\Lambda_0}^{1}) \,
 \sum _{i \in \Lambda_0} \,H_{1}^{i, \Lambda} 
\, \, p(s, u_{\Lambda_0}^{1}, \Gamma)\, \, d u_{\Lambda_0}^{1} 
\, d u_{\Lambda_0} \, ds.\]
Finally from (\ref{re.4}) we obtain the formula for the density defined by (\ref{red}):
\begin{equation}\label{re.5}
p_{A}(t, v_{A})= p^{0}_{A}(t, v_{A}) \, 
\end{equation}

\[ + \, \sum _{k=1}^{\infty} 
\, \sum _{(i_{1}, ..., i_{k}) \in \Lambda_0 ^{k}} \, 
\int _{0}^{t} \int _{{\bf R}_{+}^{\Lambda_0}} \,  
\int _{{\bf R}_{+}^{\Lambda_0}} \, \left( 
{\prod} _{z \in \Lambda_0 } g_3(u_{z}) \right)
\,p^{0}(t-s_{1}, u_{\Lambda_0}, u_{\Lambda_0}^{1}) \, \]
\[H_{1}^{i_{1}, \Lambda_0} \, \int _{0}^{s_{1}} 
\int _{{\bf R}_{+}^{\Lambda_0}}\,p^{0}(s_{1}-s_{2}, u_{\Lambda_0}^{1}, 
u_{\Lambda_0}^{2}) \ldots \]

\[ H_{1}^{i_{k}, \Lambda_0} \, \int _{0}^{s_{k-1}} 
\int _{{\bf R}_{+}^{\Lambda}}\,p^{0}(s_{k}, u_{\Lambda_0}^{k}, v_{\Lambda_0})\, 
d v_{\Lambda_0 \setminus A} \, d u_{\Lambda_0}^{k} \, 
ds_{k} \ldots du_{\Lambda_0}^{1} \, du_{\Lambda_0} \, ds_{1}, \]
Making use of the formula (\ref{ren}) and the bounds (\ref{as1}) and (\ref{p1}), 
one can prove the convergence of the series in (\ref{re.5})
and get the necessary bounds by following the 
 the formulae in Section 3 of  \cite{t1} (p.179 and further)
with only minor modifications due to (\ref{re.1}). 
Therefore 
for the sake of brevity
we  skip the rest of the straightforward part of the
 proof of Lemma \ref{L2}.

\bigskip

Consider now the embedded Markov chain $(x_i(n):=X_i^{  \Lambda_0 }(\tau _n),
i \in \Lambda_0 )$,
$n=1,2, \ldots, $ where $\tau_1, \tau _2, \ldots$ denote the 
consecutive moments when at least one of the components
of $X^{  \Lambda_0 }$ reaches zero, i.e. 
$X_j^{  \Lambda_0 }(\tau _n)=0$ for
some $j \in \Lambda_0 $. Due to ergodicity of $(X_i^{  \Lambda_0 }(t), \ 
i \in \Lambda_0)$,
there exists a limiting distribution, 
call it $F$, of 
$x(n)$ as $n \rightarrow \infty$. 
Let $x^{\infty}(n)$ be the stationary version of 
the Markov chain $x(n)$, i.e. whose initial distribution is $F$. 
 Further
 for any $i \in \Lambda_0 $ and
$ j \in D_E(i) $ let $Z_{ij}$
be  a random variable
with  distribution function
$$F_{ij}(u):=P\{x^{\infty}_i(n) <u \ \mid
\ x^{\infty}_j(n)=0 \}.$$ 
Clearly, the density function $f_{ij}$ for each
$F_{ij}$ is given by
\begin{equation}\label{S1}
f_{ij}(u) = \frac{p_{\{ij\}}(u,0)}{p_j(0)}
\end{equation}
as long as $p_j(0) >0$. Observe that $p_j^0(0)=\frac{1}{{\bf E}Y}=1$ by (\ref{1.70}).
Hence, according to Lemma \ref{L2} condition $p_j(0) >0$ is satisfied 
 at least for all small values of $w_E$.

We shall derive now an equation for $\pi_i^ { \Lambda_0, 0
}=\pi^{+, 0}(w_E )$ (see definition (\ref{f0}) and (\ref{1.5}) ).
 From (\ref{1.101}) and the above arguments on ergodicity  we obtain for all $i \in 
\Lambda_0$:
\begin{equation}\label{2.101}
 \lim _{t \rightarrow \infty }\frac{1}{t}X_i^{ 
\Lambda_0 }(t) =0=-1 + \pi^{+, 0}(w_E ) {\bf E}Y -K_E \pi^{+, 0}(w_E)
w_E
\end{equation}
\[-\sum _{j \in D_E(i) } \pi^{+, 0}(w_E)  {\bf E}( Z_{ij} -Y-
w_E \eta _{2})I\{ Z_{ij} \leq  w_E \eta _{2} \},\]
where $Z_{ij}$, $Y$ and $\eta _{2}$ are independent by their definitions.
Therefore taking into account that ${\bf E}Y=1$,
we obtain from (\ref{2.101}) the following equation:
\[ \pi^{+, 0}(w_E) \left(1- K_E w_E + \sum _{j \in D_E(i) } 
{\bf P} \{ Z_{ij} \leq  w_E \eta _{2} \}
\right.\]
\begin{equation}\label{1.10}
\left. + \sum _{j \in D_E(i) }\  {\bf E}(w_E \eta _{2}-Z_{ij})
I\{ Z_{ij} \leq  w_E \eta _{2} \} 
\right) \ = \ 1.
\end{equation}

\begin{lem} \label{L1} There exists a positive
constant $C_1$ independent of $N$  such that for any $i \in \Lambda_0 $
and for any
$ j \in D_E(i) $
\begin{equation}\label{S3}
| {\bf P}\{Z_{ij} \leq  w_E\eta _{2}\}  \ - \ w_E| \leq w_E ^{2}C_1
\end{equation}
and
\begin{equation}\label{S4}
{\bf E}(w_E \eta _{2}-Z_{ij}) I\{ Z_{ij} \leq  w_E \eta _{2} \} \leq w_E ^{2}C_1
\end{equation}
for all $0 \leq w_E \leq c$ (with $c$  defined in Lemma \ref{L2}).
\end{lem}

\noindent
{\bf Proof of Lemma \ref{L1}.}
In the case $w_E=0$ let $Z_{ij}=Z_{ij}^0$. Then the density
 of $Z_{ij}^0$ (call it correspondingly  $f_{ij}^0(u)$ ) is given by
\begin{equation}\label{S2}
f_{ij}^0(u)=p^0(u)=\int_u^{\infty}g_0(y)dy
\end{equation}
due to (\ref{S1}) and (\ref{1.70}).

First let us obtain the following bound as a corollary of Lemma \ref{L2}:
\begin{equation}\label{1.130}
| f_{ij}(u)- f_{ij}^0(u) | \leq w_E C_2, \ \ u \in {\bf R}_+,
\end{equation}
where $C_2$ is some
positive constant independent of $N$. Notice, that from (\ref{2.L2}) and  (\ref{1.7})
we immediately derive:
\begin{equation}\label{1.13}
| p_{A}(u_{A})- p_A^0(u_{A}) | \leq w_E C_3 \exp \{ -\beta \sum_{i \in A} u_i\}, 
\end{equation}
for any $u_{A}\in {\bf R}_+^{A}$, 
where $C_3 $ is some positive constant independent of $N$. 
Then the bound (\ref{1.130}) readily follows from (\ref{1.13})
and the continuity of the transformation (\ref{S1}).

Let us  prove now  (\ref{S3}). Consider 
\begin{equation}\label{*1}
 {\bf P} \{Z_{ij} \leq  w_E\eta _{2}
\} =\int _0^{\infty}g_2(y) \ \int _0^{w_Ey} f_{ij}^0(u)du \ dy 
\end{equation}
\[+ \int _0^{\infty}g_2(y) \ \int _0^{w_Ey}(f_{ij}(u) - f_{ij}^0(u)) \ du \ dy.\]
Substituting formula (\ref{S2}) into the first integral in (\ref{*1}) we derive
after simple calculations:
\begin{equation}\label{*2}
\int _0^{\infty}g_2(y) \ \int _0^{w_Ey} f_{ij}^0(u)du \ dy 
\end{equation}
\[ = \int _0^{\infty} {\bf P} \{\eta _2 \geq y \} w_E \ (1- {\bf P} \{Y \leq 
w_Ey \})\ dy \]
\[ = w_E - w_E \int _0^{\infty} {\bf P} \{\eta _2 \geq y \} {\bf P} \{Y \leq 
w_Ey \}\ dy, \]
where we used the condition ${\bf E}\eta _2=1$. Taking into account
(\ref{as1}) it is easy to see that
\begin{equation}\label{*3}
w_E \int _0^{\infty} {\bf P} \{\eta _2 \geq y \} {\bf P} \{Y \leq 
w_Ey \}\ dy  \leq C'w_E^2
\end{equation}
for some positive constant $ C'$.
Next, making use of (\ref{1.130})  we readily derive 
\begin{equation}\label{*4}
|\int _0^{\infty}g_2(y) \ \int _0^{w_Ey}(f_{ij}(u) - f_{ij}^0(u)) \ du \ dy|
\leq    w_E^{2} C_2.
\end{equation}
Combining the bounds (\ref{*4}) and  (\ref{*3}) together with 
equations  (\ref{*2})
and (\ref{*1}) gives us (\ref{S3}).

To prove  (\ref{S4}) let us consider the following decomposition:
\[ {\bf E}(w_E \eta _{2}-Z_{ij}) I\{ Z_{ij} \leq  w_E \eta _{2} \} \]
\[ = \int _0^{\infty}g_2(y) \left( \int _0^{w_Ey}(w_Ey-u) f^0_{ij}(u) \, du \right)
\ dy\]
\[ + \int _0^{\infty}g_2(y) \left( \int _0^{w_Ey}(w_Ey-u) (f_{ij}(u) - f_{ij}^0(u))
\, du \right)
\ dy.\]
Using again (\ref{as1}) and (\ref{1.130}) we immediately derive from here:
\begin{equation}\label{*5}
{\bf E}(w_E \eta _{2}-Z_{ij}) I\{ Z_{ij} \leq  w_E \eta _{2} \} 
\leq  w_E^2 C_3 +w_E^3 C_4,
\end{equation}
 where $C_3 $ and $C_4$ are some positive constants, which implies (\ref{S4}).
The lemma is proved.

Lemma \ref{L1} allows us to derive from  equation (\ref{1.10}) that
\begin{equation}\label{*6}
|\pi^{+, 0}(w_E)-1| \leq  w_E^{2}C_5
\end{equation}
for some positive constant $C_5$.

Finally let us consider $\pi^{+, e}_{ij}(w_E)$ for $i \in \Lambda _0, 
j \in D_E(i).$
According to  definition (\ref{fe}) we derive similarly  (\ref{2.101})
\begin{equation}\label{*7}
\pi _{ij}^{W, e}= \lim _{t \rightarrow \infty}
\frac{1}{t} \# \{0<t_n<t: \  X_j^{ W}(t_n)=0 \mbox{ and } X_i^{ W}(t_n)
- \theta^{(n)}_{ji} \leq 0\}
\end{equation}
\[ =\lim  _{t \rightarrow \infty}
\frac{1}{t} \sum _{0<t_n <t: \ X_j^{ \Lambda_0 }( t_n) = 0 } 
 I\{ X_i^{\Lambda_0 }(t_n ) \leq w_E \eta _{2}^{(n)} \}\]
\[ = \pi^{+, 0}(w_E) \ {\bf P} \{Z_{ij} \leq  w_E\eta _{2}\}.\]
Substituting (\ref{*7}) into  definition  (\ref{1.5}) we obtain
\begin{equation}\label{*8}
\pi^+(w_E)=\pi^{+, 0}(w_E) (1+ \sum _{j \in D_E(i) } 
{\bf P}\{Z_{ij} \leq  w_E\eta _{2}\}),
\end{equation}
which together with  (\ref{*6}) and (\ref{S3}) gives us
the following bound:
\begin{equation}\label{1.14}
|\pi^+(w_E) - (1 + K_Ew_E)| \leq w_E^{2} C_6
\end{equation}
for all $w_E \leq C_0 $, where $C_6$ and $C_0$ are some positive constants.

Substituting the bound (\ref{1.14}) into (\ref{CR}) we get (\ref{21.4}). This
completes the proof of Theorem \ref{T2}.

\section{Fully connected finite networks.}
\setcounter{equation}{0}

\subsection{The model.}
Let $\Lambda
=\{1, \ldots , N\}$.
Consider an $N$-neuron network 
with inhibitory connections only, i.e.
$X(t)=X^{inh}(t)$ as defined in (\ref{ID}). Here we assume that
$D(i)=\Lambda \setminus \{i\}$ 
for any $i \in \Lambda $, which means that our network is
{\it fully connected}. 

Let $N=2pk$, where $k,p \in Z_+$ are fixed arbitrary.
Suppose a network consists of $2p$ symmetric subnets (blocks), 
determined by the values of the
connection constants ${\bf E}{\theta _{ij}}$ as follows.
We divide set $\Lambda =\{1, ... , N\}$ into $2p$ 
subsets $W_i$, $i=1,..., 2p$, so that 
\begin{equation}\label{c0}
|W_i|=k \mbox{ and } \cup _{i=1}^{2p} W_i =\Lambda .
\end{equation} 
Notice that we keep $k>1$ fixed but arbitrary. This is simply a size of one
block or unit of our network. The case $k=1$ is  trivial.

Let us divide also the set $\{1, ... , 2p\}$ into a set of $p$ non-intersecting 
pairs, i.e., 
\begin{equation}\label{c1}
\{1, ... , 2p\}= \cup_{n=1}^{p} C_n= \cup_{n=1}^{p} \{C_n^1,C_n^2\},
\end{equation}
where
\[|C_n|=2, \mbox{ and } C_i \cap C_j=\emptyset. \]   
Then we set
\begin{equation}\label{3.1}
\begin{array}{ll}
{\bf E} Y_x  & =a>0,\\ 
{\bf E} \theta _{xy} & = -c_{ij} <0,  \mbox{ if }  x \in W_i, \, \, y \in
W_j, \end{array}
\end{equation}
for all $x,y \in \Lambda $ and $1 \leq i,j \leq 2p$, where
\begin{equation}\label{c3}
c_{ij}= \left\{ 
\begin{array}{ll}
c, & \mbox{ if } \{i,j \} = C_n \mbox{ for some }n, \\
b, & \mbox{ otherwise, }
\end{array}
\right.
\end{equation}
for some arbitrary but fixed constants
\begin{equation}\label{c}
0<b<a<c.
\end{equation}

Thus our network consists of $2p$ connected blocks of interacting identical
neurons. The connections between the neurons of different blocks are
different in general from those between the
neurons of the same block. 
In the following Theorem \ref{p2} below, we describe all the possible
 traps for our network, i.e. all the possible limiting states in the transient case.
 Notice that the
total number of the patterns is  $|T|=2^p=2^{\frac{N}{2k}}$,
 which is exponentially large 
with respect to $N$.
This description will allow us to find a method to reconstruct or "memorize",
in some particular cases, a given system of traps
 (see Corollary \ref{cor1} below).
More precisely, we will show how to choose the appropriate connection constants.

\subsection{Description of patterns.}

\begin{theo}\label{p2}
Let $N=2pk$, where $k>1,p>1$, and decompositions 
(\ref{c0}) and (\ref{c1}) be fixed but arbitrary.
Call ${\cal L}$ the set of all the decompositions of $\{1, ... , 2p\}$ 
into a pair of subsets such that
\begin{equation}\label{c2}
{\cal L}:= \{(L_1, L_2): \  L_1 \cup L_2 = \{1, ... , 2p\}, \
L_1 \cap L_2 = \emptyset, 
\end{equation}
\[  L_i \cap C_n =1, \ \ i=1,2, \ n=1, ..., p\}.\]
Then a set $A \subset \Lambda$ is a trap for the network with 
the connection constants 
(\ref{3.1})-(\ref{c}), if and only if
\begin{equation}\label{c4}
A=\{x \in \Lambda: \ x \in \cup_{i \in B} W_i\},
\end{equation}
for some $B \subset \Lambda$ such that
\begin{equation}\label{c5}
(B, \{1, ... , 2p\} \setminus B) \in {\cal L}.
\end{equation}
\end{theo}

{\bf Proof.}
First we will show that any set $A$ satisfying conditions  (\ref{c4}) 
and (\ref{c5}) is a trap for the above defined system. 
Indeed,  consider a restriction $X^{\Lambda \setminus A}(t)$.
In this case $(X_i^{\Lambda \setminus A}(t), \ i \in \Lambda \setminus A)$
 is a completely 
connected system, and for any $x,y \in \Lambda \setminus A$, $y \neq x$
\begin{equation}\label{c6}
- {\bf E} \theta _{xy}=b< a={\bf E} Y_x,
\end{equation}
which implies ergodicity of $(X_i^{\Lambda \setminus A}(t), 
\ i \in \Lambda \setminus A)$ due to Proposition 2.2 in
\cite{kmr}. 

Assume, $(X_i^W(t), \ i \in W)$ is ergodic for some $W \subset \Lambda$.
Observe that in the case of a fully connected network with only
negative connections  definition (\ref{2.2}) 
becomes
\begin{equation}\label{2.2*}
\pi _i^W=\pi _i^{W,0}, \ i \in W.
\end{equation}
Furthermore it follows by ergodicity from (\ref{ID}) and definition (\ref{f0})
that  
\begin{equation}\label{2.3*}
\pi _i^W=\left( {\bf E}Y_i + \sum_{j \in W \setminus \{i\}}{\bf E}|\theta _{ji}|
 \right)^{-1}, 
\ i \in W.
\end{equation}
Correspondingly, we derive from (\ref{2.1})
\begin{equation}\label{2.1*}
v_j^W =-1 + \sum_{i \in W}{\bf E}|\theta _{ij}| \pi _i^W, \ j \in \Lambda \setminus W.
\end{equation}

Notice that definitions (\ref{c0})-(\ref{c1}) imply $|B|=p$
for the set $B$ satisfying conditions (\ref{c4})-(\ref{c5}).
Then for any $x \in \Lambda \setminus A$ we derive from 
(\ref{2.3*}) and condition (\ref{c3}) (see also \cite{kmr}), that
\begin{equation}\label{c7}
\pi _{x}^{\Lambda  \setminus A}=\frac{1}{a+(|\Lambda \setminus A|-1)b}
=\frac{1}{a+(pk-1)b}.
\end{equation}
Substituting (\ref{c7}) into definition (\ref{2.1*}) and taking into account
(\ref{3.1})-(\ref{c3}), we obtain for any $y \in A$
\begin{equation}\label{c8}
v_y^{\Lambda \setminus A} := -1  + \frac{ck+(p-1)bk}{a+(pk-1)b},
\end{equation}
which together with (\ref{c})
implies
\[ v_y^{\Lambda \setminus A}>0.\]
Hence, we conclude that $A$ is a trap.

Next we will show that any trap satisfies (\ref{c4})-(\ref{c5}).
Clearly, any subset $A \subseteq \Lambda$ can be represented as
\begin{equation}\label{c9}
A=\{x \in \Lambda: \ x \in \cup_{i \in B} A_i\},
\end{equation}
where $B \subseteq \{1, ... , 2p\}$ and $A_i \subseteq W_i$ for any $i \in B$.
Suppose a subset $A$ does not satisfy conditions (\ref{c4})-(\ref{c5}),
which happens if and only if
at least one of the following situations (I) or (II) takes place:\\
\noindent
(I) the set $B$ does not satisfy (\ref{c5}), which means that
\begin{equation}\label{c10}
B=\{C_n, \ n \in I \} \cup \{C_n^{l}, \ n \in J\}
\end{equation}
and
\begin{equation}\label{c11}
{\bar B}:=\{1, \ldots , 2p\} \setminus B = \{C_n, \ n \in I' \} \cup
\{C_n^{l'}, \ n \in J\},
\end{equation}
where $\{l,l'\}=\{1,2\}$ and sets $I$ and $I'$ are such that
\begin{equation}\label{c12}
I \cup I' \neq \emptyset, \ I \cap I'=\emptyset, \ \mbox{ and } I \cup I'
\cup J =\{1, \ldots , 2p\};
\end{equation}

\noindent
(II) $W_i \setminus A_i \neq \emptyset$ at least for some $i \in B$, i.e. 
(\ref{c4}) is not satisfied.

Suppose situation (I) takes place. Without loss of generality let $l=1$ and $l'=2$. 
Consider
\begin{equation}\label{c13}
{\bar A}:= \Lambda \setminus A =
\{ \cup_{i \in {\bar B}} W_i\} \cup \{ \cup_{i \in B} W_i \setminus A_i \}.
\end{equation}
Define the subset $B_1 \subseteq B$ so that 
\begin{equation}\label{B1}
W_i \setminus A_i \neq \emptyset \mbox{ iff } 
i \in B_1,
\end{equation}
and thus 
\[{\bar A}= \{ \cup_{i \in {\bar B}} W_i\} \cup \{ \cup_{i \in B_1} W_i
\setminus A_i \}. \]
 Notice that $B_1$ can be empty. 
Further define $I_1 $ so that
\begin{equation}\label{c14}
C_n \in {\bar B} \cup B_1 \ \mbox{ iff } n \in I_1.
\end{equation}

(a) Suppose $I_1 \neq \emptyset$. Let 
\begin{equation}\label{c15}
B_0:=\{C_n^1, \ n \in I_1\} \cup \{{\bar B} \cup B_1 \setminus \{C_n, \ n
\in I_1\} \}.
\end{equation}
Denote
\[ {\bar A}^{B_0}=\{ \cup_{i \in {\bar B}\cap B_0} W_i\} \cup \{ \cup_{i \in
B_1 \cap B_0} W_i \setminus A_i \}. \]
Then the subsystem  $(X_i^{{\bar A}^{B_0}}(t), \ i \in {\bar A}^{B_0})$ 
is ergodic, since it is
completely connected, and 
 (\ref{c6}) holds for any $x, y \in {\bar A}^{B_0}$, $x \neq y$. 
Also, we can compute as in (\ref{c7})
\begin{equation}\label{c161}
\pi _x^{{\bar A}_{B_0}}=  \frac{1}{a+(|{\bar A}^{B_0}|-1)b},
\end{equation}
Taking into account that $|I_1| \geq 1$, we easily derive from (\ref{c161})
the following upper bound for the components of the second vector field 
$v^{{\bar A}^{B_0}}$:
\begin{equation}\label{c16}
v_z^{{\bar A}^{B_0}} \geq  -1  + \frac{c+(|{\bar A}^{B_0}|-1)b}{a+(|{\bar
A}^{B_0}|-1)b}>0
\end{equation}
 for any $ z \in {\bar A}\setminus {\bar A}^{B_0} $ due to condition
(\ref{c}). Hence, $(X_i^{\bar A}(t), \ i \in {\bar A})$ is 
transient according to Theorem A. This contradicts our assumption that
$A$ is a trap.

(b)  Suppose now that $I_1 = \emptyset$. Then it follows from (\ref{c14}) 
and (\ref{c11})
that we also have
\begin{equation}\label{c19}
I' = \emptyset. 
\end{equation}
  This together with our assumption (\ref{c12}) implies that
$I \neq \emptyset $. Thus in this case we have
\begin{equation}\label{c17}
A=\{x \in \Lambda: \ x \in \cup_{i \in B} A_i\},
\end{equation}
where $B \ni \{C_n, \ n \in I \}$, and
\begin{equation}\label{c18}
{\bar A}=
\{ \cup_{i \in {\bar B}} W_i\} \cup \{ \cup_{i \in B_1} W_i \setminus A_i \}, 
\end{equation}
where $\{ {\bar B} \cup B_1 \} \cap C_n \leq 1$ for any $n$, according to
(\ref{c14}) 
and (\ref{c19}). The latter 
implies that  (\ref{c6}) holds 
for any $y \in {\bar A}$ which in turn implies ergodicity of $(X_i^{\bar A}(t), 
\ i \in {\bar A} )$.
Furthermore, we can find $\pi _x ^{\bar A}$, analogously to (\ref{c7}), namely:
\begin{equation}\label{c20}
\pi _{x}^{\bar A}
=\frac{1}{a+(|{\bar A}|-1)b}, \ \ x \in {\bar A}.
\end{equation}
Let us compute now the $x$th component of the second vector field $v_x^{\bar
A}$ for $x \in A_i$
with $i \in I$. According to (\ref{c20}) and assumption (\ref{c3}) we get
\begin{equation}\label{c21}
v_x^{\bar A} := -1  + \frac{|{\bar A}|b}{a+(|{\bar A}|-1)b}<0.
\end{equation}
The latter contradicts our 
assumption, that $A$ is a trap.

Hence we conclude that if $A$ is a trap then necessarily condition
(\ref{c5}) holds. Assume now situation (II). More precisely, suppose that
$A$ satisfies (\ref{c5}) but does not satisfy (\ref{c4}),
i.e.
\begin{equation}\label{c22}
A=\{x \in \Lambda: \ x \in \cup_{i \in B} A_i\},
\end{equation}
for some $B \subset \Lambda$ such that
(\ref{c5}) holds, while
 $W_i \setminus A_i \neq \emptyset$   for  $i \in B_1 \subseteq B$ for some
nonempty $B_1$ (see definition (\ref{B1}) ).
The latter implies
\[
{\bar A}=
\{ \cup_{i \in {\bar B}} W_i\} \cup \{ \cup_{i \in B_1} W_i \setminus A_i \},
\]
where according to the assumption (\ref{c5})
\begin{equation}\label{c24}
| B_1 \cap C_n |=1
\end{equation}
for at least one $C_n$. The latter 
implies $C_n \in  {\bar B} \cup B_1$, i.e. the set $I_1$ defined in (\ref{c14})
is non-empty. But as we have seen, this situation contradicts  the
assumption that $A$ is a trap. This finishes our argument that conditions
(\ref{c4}) and
(\ref{c5}) are necessary for a set $A$ to be a trap. This completes the
proof of Theorem \ref{p2}.

\subsection{Learning rule for the almost symmetric case}

We will show here that the Hebbian rule (see \cite{H}) of learning
 patterns which successfully
 works
in the case of Hopfield networks (\cite{H}), is applicable for our network
at least for the following particular case.

For any $A \subset \Lambda $ define
configuration
 $\xi (A)=(\xi _i (A), \ i \in \Lambda )\in \{-1, +1
\}^{\Lambda }$ such that
\begin{equation}\label{TD}
\xi _i (A) = \left\{
\begin{array}{ll}
+1, & \mbox{ if } i \in  A, \\
-1, & \mbox{ otherwise. }
\end{array}
\right.
\end{equation}
We shall also call a configuration
 $\xi$ a trap for $X(t)$ if and only if $\xi=\xi (A)$, where
$A$ is a trap.

Suppose that we are given $2^p$, $p =\frac{N}{2k}$,  binary vectors (images)
 $\xi _1, \ldots , \xi _{2^p}$.
We shall find the 
connection constants ${\bf E}\theta _{xy}$ such that $X(t)$ 
satisfying (\ref{3.1}) with these
 parameters, possesses a system of traps consisting exactly of the given  
$2^p$ vectors. 

Notice the difference between this 
task and the problem of stability of patterns for  Hopfield neural model 
 (see for example, \cite{AB} for a recent account on the
relevant results). Recall, that the capacity of Hopfield network is
determined by the number of given {\it i.i.d.} patterns, 
which are stable with respect to the
dynamics of the system. This means that starting from an arbitrary initial state 
the system should converge with a large probability to one of the given patterns,
 which is the
 closest to the initial state. It was conjectured that the
number of such patterns for Hopfield network of $N$ neurons is at most a fraction
of $N$. Here we construct a network  which posesses  given exponentially large 
(with respect to $N$) number of the limiting patterns and only them. These patterns are
stable in a trivial sense, i.e. if 
the initial state of the system is one of the 
given traps (patterns), then the system stays at this trap forever.
However, we do not predict which state (out of $p$ possible) the
system converges to, starting from 
an arbitrary initial condition. The problem of 
determination of the basins of 
attraction of the limiting patterns for the hourglass model
will be a subject of a separate study.

\begin{cor} \label{cor1} Suppose the collection of $N$-dim vectors 
$\{\xi ^{\mu}, \mu = 1, \ldots, M\}$ where $M=2^p$,
 and $\xi _x^{\mu} \in \{-1,1\} $ for any $\mu$ and $1
\leq x \leq N$, has the following properties:
\medskip

1. $\xi _x^{\mu}\xi _y^{\mu} =1$ for any $\mu$ if $x,y \in  W_n$ for some 
$n $ $\in \{1,...,2p\}$,
\medskip

2. $\sum_{x=1}^{N}\xi _x^{\mu}=0$ for all $\mu$,
\medskip

3. for any $n \in \{1,...,2p\}$ there exists unique $ l=l(n) \in
\{1,...,2p\} \setminus\{n\}$ such 
that
  $\xi _x^{\mu}\xi _y^{\mu} =-1$ for any $\mu$ whenever $x \in  W_n$ and $y
\in  W_l$. 
\medskip

Then the $N$-neuron system with $a_i=a$ and the connection constants
\begin{equation}\label{lp1}
 {\bf E}\theta _{xy} = \left\{
\begin{array}{ll}
b(x,y), & \mbox{ if } b(x,y)= \min_{(x',y')} b(x',y')\\
\max_{(x',y')} b(x',y'), & \mbox{ otherwise, }
\end{array}
\right.
\end{equation}
where
\begin{equation}\label{lp3}
 b(x,y) := Aa \frac{1}{M}\, \sum _{\mu=1}^{M}\xi _x^{\mu}\xi _y^{\mu} \, -
\, Ba, 
\end{equation}
with the constants $A$ and $B$ satisfying the conditions
\begin{equation}\label{lp2}
\begin{array}{rll}
0 < & B-A & < 1,\\
1 < & B+A, & 
\end{array}
\end{equation}
has a system of traps, which is $\{\xi ^{\mu}$, $\mu = 1, \ldots, 2^p \}$.
\end{cor}

{\bf Proof.} Indeed, 
 having conditions {\it 1-3} of the corollary satisfied, we
 derive from (\ref{lp3}) that
\begin{equation}\label{lp4}
  b(x,y) = 
\left\{
\begin{array}{ll}
-(A+B)a, & \mbox{ if } x \in W_n, y \in W_{l(n)}, \\
-(B-A)a, & \mbox{ if } x , y \in W_n,\\
\end{array}
\right.
\end{equation}
and
\begin{equation}\label{lp5}
 -(A+B)a < b(x,y) < -(B-A)a \ \  \mbox{ otherwise.}
\end{equation}
Substituting (\ref{lp4}) snd (\ref{lp5}) into 
(\ref{lp1}), and taking into account condition (\ref{lp2}), we obtain
\begin{equation}\label{lp7}
\begin{array}{ll}
-{\bf E}\theta _{xy} =(A+B)a >a, & \mbox{ if } x \in W_n, y \in W_{l(n)},\\
-{\bf E}\theta _{xy} =(B-A)a <a, & \mbox{ otherwise, }
\end{array}
\end{equation}
which shows that conditions (\ref{c3})-(\ref{c}) are satisfied. Hence we can use
 Theorem \ref{p2} to conclude that $\{\xi ^{\mu}$, 
$\mu = 1, \ldots, 2^p \}$ is a system of traps 
for the defined network.

\medskip
{\bf Acknowledgements} We thank the referee for useful remarks.

\end{document}